\documentclass[12pt,reqno]{amsart}
\usepackage{amssymb}
\usepackage{amsmath,amssymb,color}
\newtheorem{thm}{Theorem}[section]
\newtheorem{lem}[thm]{Lemma}
\newtheorem{cor}[thm]{Corollary}

\numberwithin{equation}{section}
\newcommand{\UUUUU}{\R^n \setminus \ooo{U}}

\newcommand{\ep}{\varepsilon}

\newcommand{\va}{\varphi}
\newcommand{\ppp}{\partial}

\newcommand{\intBB}{\int_{B_{2R}}}
\newcommand{\ddd}{\mbox{div}\thinspace}

\newcommand{\weight}{e^{2s\va}}

\newcommand{\sumij}{\sum_{i,j=1}^n}
\newcommand{\sumk}{\sum_{k=1}^n}

\newcommand{\R}{\mathbb{R}}

\newcommand{\N}{\mathbb{N}}

\newcommand{\www}{\widetilde}

\newcommand{\ooo}{\overline}

\newcommand{\RRRMN}{(\R^m\setminus \{ \vert y\vert \le 1\})\times \R^{n-m}}

\newcommand{\YYYY}{\vert y\vert}
\newcommand{\ZZZZ}{\vert z\vert}
\allowdisplaybreaks
\begin{document}
\title
[]
{Uniqueness of solutions to an elliptic inequality with 
rapid decay at infinity}

% with 
%time order $\alpha \in (1,2)$.
%}

\pagestyle{myheadings}
\markboth{F. G\"olgeleyen, O.Y. Imanuvilov and M. Yamamoto}{F. G\"olgeleyen, 
O.Y. Imanuvilov and M. Yamamoto}

\author{
$^1$ F. G\"olgeleyen, $^2$ O.Y. Imanuvilov and $^{3}$ M.~Yamamoto
}

\thanks{
$^1$ Department of Mathematics, Faculty of Science, 
Zonguldak B\"ulent Ecevit University, Zonguldak 67100 T\"urkiye
e-mail: {\tt f.golgeleyen@beun.edu.tr}
\\
$^2$ Department of Mathematics, Colorado State University, 101 Weber Building, 
Fort Collins CO 80523-1874, USA 
e-mail: {\tt oleg@math.colostate.edu}
\\
$^3$ 
Department of Mathematics, Faculty of Science, 
Zonguldak B\"ulent Ecevit University, Zonguldak 67100 T\"urkiye\\
Graduate School of Mathematical Sciences, The University
of Tokyo, Komaba, Meguro, Tokyo 153-8914, Japan 
e-mail: {\tt myama@ms.u-tokyo.ac.jp}
}

\date{}
%\begin{document}
\maketitle

\baselineskip 18pt

\begin{abstract}
We consider an elliptic differential inequality: 
$\vert \Delta u(x) \vert \le C_0(\YYYY^{-\gamma}\vert u(x)\vert 
+ \YYYY^{-\theta}\vert \nabla u(x)\vert)$ in an 
exterior domain $\R^n \setminus \ooo{U}$, where $U$ is a simply connected 
bounded domain $U$,
$x := (y,z) \in \R^n$ with $y \in \R^m$ and $z\in \R^{n-m}$
for given $m\in \{ 1, ..., n\}$, and 
$\gamma, \theta \in \R$ are constants.
We assume that $u(x)$ decays with 
exponential rate in the $y$-coordinates and polynomial rate in 
the $z$-coordinates as $\vert x\vert \to \infty$. 
We prove that if decay rates of $u$ satisfy certain conditions related 
to the constants 
$\gamma, \theta \in \R$, then $u\equiv 0$ in $\UUUUU$.  
The key is a Carleman estimate with typical cut-off arguments.
%\\
%{\bf Key words.}  
%time-fractional diffusion equation, backward problem,
%well-posedness, divided data at muiltiple moments 
%\\
%{\bf AMS subject classifications.}
%26A33, 34A08, 34A12
\end{abstract}

\section{Introduction}

We start with an elliptic equation in $\UUUUU$, where 
$U \subset \R^n$ is a simply connected bounded domain:
\begin{equation}
\Delta u(x) + A(x)\cdot \nabla u(x) - V(x)u(x) = 0 
\quad \mbox{for $x \in \UUUUU$}.             \label{1.1}  
\end{equation}
Henceforth we set $H^2_{loc}(\UUUUU):= \{ u;\, u\vert_M 
\in H^2(M)\,\, \mbox{for any bounded domain $M \subset \UUUUU$}\}$.

We are concerned with decay rates of the solution $u\in 
H^2_{loc}(\UUUUU)$ as $\vert x\vert \to \infty$ 
which imply that $u \equiv 0$ in $\UUUUU$.

In a special case where $A=0$ and $V=0$ in $\R^n$, that is,
$u \in H^2_{loc}(\R^n)$ is harmonic in $\R^n$, 
the classical Liouville theorem reads 
\begin{equation}
\lim_{\vert x\vert \to\infty} \vert u(x)\vert = 0  \label{1.2}
\end{equation}
implies that $u=0$ in $\R^n$.  We understand that \eqref{1.2} is the weakest 
decay condition, but \eqref{1.2} does not imply $u=0$ in $\R^n$ in general 
with $A\not\equiv 0$ or $V\not\equiv 0$.

For an equation
\begin{equation*}
\Delta u(x) - V(x)u(x) = 0 \quad \mbox{for $x \in \R^n$}
\end{equation*}
with $V\in L^{\infty}(\R^n)$, Landis posed a conjecture that if 
there exist constants $C_1, C_2 > 0$ and $\ep>0$ such that
\begin{equation*}
\vert u(x)\vert \le C_1\exp(-C_2\vert x\vert^{1+\ep}) \quad 
\mbox{for $x\in \R^n$},
\end{equation*}
then $u=0$ in $\R^n$ (Kondratiev and Landis \cite{KL}).  
Meshkov \cite{M} demonstrated that Landis' conjecture does not hold 
by constructing a complex-valued solution $u\not\equiv 0$ satisfying 
\begin{equation*}
\vert u(x)\vert \le C_1\exp\left(-C_2\vert x\vert^{\frac{4}{3}}\right) \quad 
\mbox{for $x\in \UUUUU$}
\end{equation*}
and that if there exist constants $C_1, C_2 > 0$ and $\ep>0$ such that 
a solution $u$ to \eqref{1.2} satisfies 
\begin{equation}
\vert u(x)\vert \le C_1\exp\left(-C_2\vert x\vert^{\frac{4}{3}+\ep}\right) 
\quad \mbox{for $x\in \UUUUU$},         \label{1.3} 
\end{equation}
then $u = 0$ in $\R^n$.  In other words, $\frac{4}{3}$ is the 
"critical decay rate for uniqueness" of solution $u$.

After paper of  Meshkov \cite{M}, there have been many important works 
concerning Landis'
conjecture.  As for the two dimensional case 
$n=2$, we can refer to Davey \cite{D2}, 
Davey, Kenig and Wang \cite{DKW1}, \cite{DKW2}, Davey and Wang \cite{DW},
Kenig, Silvestre and Wang \cite{KSW}, Kenig and Wang \cite{KW},
Le Balc'H and de Souza \cite{LBdS}, Logunov, Malinnikova, Nadirashvili and 
Nazarov \cite{LMNN}.
As for $n\ge 3$, see Davey \cite{D1}, \cite{D3},
Davey and Zhu \cite{DZ}.
In particular, for the equation $\Delta u + A(x)\cdot \nabla u - V(x)u = 0$ 
with $n\ge 3$, Theorems 2-4 in Davey and Zhu \cite{DZ} give
the critical decay rate for uniqueness according to Lebesgue spaces to 
which $A$ and $V$ belong, while Davey \cite{D1} can provide the 
critical decay rates for uniqueness in terms of the polynomial decay
rates of $A(x)$ and $V(x)$ as $\vert x\vert \to \infty$.
See also Rossi \cite{R}.
Here we do not intend to create any comprehensive references. 

Moreover Landis' conjecture is related to the strong unique continuation
property, for which we refer, for instance, to 
Alessandrini \cite{A}, Garofalo and Lin \cite{GL}, Kenig \cite{K}.

Now we start to formulate our problem and state the main results.
Throughout the article, let functions under consideration be real-valued.
We can modify the arguments for the case of complex-valued functions.

For arbitrarily fixed $m\in \{1,2,..., n\}$, we set 
$$
x = (x_1,.., x_m,x_{m+1},..., x_n) \in \R^n, \quad
y:= (x_1,..., x_m) \in \R^m, \quad z:= (x_{m+1}, ..., x_n)
\in \R^{n-m}.
$$
We write $\vert y\vert := (x_1^2+ \cdots + x_m^2)^{\frac{1}{2}}$
and $\vert z\vert := (x_{m+1}^2+ \cdots + x_n^2)^{\frac{1}{2}}$, and 
$\ooo{U}$ denotes the closure of a set $U \subset \R^n$.
We denote $\ppp_k = \frac{\ppp}{\ppp x_k}$ for $1\le k \le n$.

In this article, we consider elliptic differential inequalities,
not elliptic equations, because our method does not need to 
assume elliptic equations.  More precisely, 
we arbitrarily choose a bounded domain $U \subset \R^n$.
We assume   
\begin{equation}
\vert \Delta u(x) \vert \le C_0(\YYYY^{-\gamma} \vert u(x)\vert 
+ \YYYY^{-\theta} \vert \nabla u(x) \vert), \quad x \in \UUUUU 
                                          \label{1.4}  
\end{equation}
with some constants $C_0 > 0$ and $\gamma, \theta \in \R$.

Our main purpose is to provide values of critical decay rates 
of $u$ which imply that $u=0$ in $\UUUUU$, in terms of 
the constants $\gamma, \theta$, and moreover we are interested in 
reducing the exponential decay $y$-coordinates.
\\

We define a constant $\alpha>0$ by  
\begin{equation}
\alpha := \max\left\{ \frac{2}{3} - \frac{\gamma}{3},\,\, 
1 - \theta, \, \, \frac{1}{2}+\ep \right\}         \label{1.5}  
\end{equation}
with arbitrarily chosen $\ep > 0$ and the constants $\gamma,\theta$  
in \eqref{1.4}.

Let $g(r) \ge 0$ for $r > 0$, satisfy  
\begin{equation}  
\int^{\infty}_0 g(r)^2r^{n-m-1} dr < \infty \quad
\mbox{if $m\in \{1, ..., n-1\}$.}                   \label{1.6}       
\end{equation}
For $m=n$, we do not need to introduce $g(r)$.
\\
{\bf Example of $g$}. If $g \in C[0,\infty)$ satisfies 
$g(r) \sim r^{\frac{m-n}{2}-\delta}$ as $r \to \infty$ with some
constant $\delta>0$, then (1.6) holds. 

We have
\\
\begin{thm}\label{T1}
{\it
We assume that $u\in H^2_{loc}(\UUUUU)$ satisfies \eqref{1.4} and 
there exist constants $C>0$, $C_1>0$ and $\beta > \alpha$ such that 
\begin{equation} 
\vert u(x)\vert \le Cg(\vert z\vert)\exp( -C_1\vert y\vert^{2\beta}),
\quad x \in \UUUUU. 
                              \label{1.7} 
\end{equation}
Then $u=0$ in $\UUUUU$.
}
\end{thm}

\begin{cor}\label{C1}
{\it We assume that $u\in H^2_{loc}(\UUUUU)$ satisfies \eqref{1.4}.
Let $m=n$ and $\beta > \alpha$, and there exist constants $C, C_1 > 0$ such 
that 
\begin{equation}
\vert u(x)\vert \le C\exp( -C_1\vert x\vert^{2\beta}),
\quad x \in \UUUUU.             \label{1.8}             
\end{equation}
Then $u=0$ in $\UUUUU$.
}
\end{cor}

\begin{cor}\label{C10}
{\it 
We assume that $u\in H^2_{loc}(\UUUUU)$ satisfies 
$\vert \Delta u(x)\vert \le C_0\vert u(x)\vert)$ for $x \in \UUUUU$.  If 
$$ 
\vert u(x)\vert \le Cg(\vert z\rvert)\exp( -C_1\vert y\vert
^{\frac{4}{3}+\ep}), \quad x \in \UUUUU
$$
with some constants $\ep > 0$ and $C, C_1 > 0$, then $u=0$ in $\UUUUU$.
}
\end{cor}

In Corollary 1.3, as is seen by the proof of Theorem 1.1, since 
we can choose the constant $\theta > 0$ in (1.4) arbitrarily,  
we can satisfy (1.5) with $\alpha = \frac{2}{3}$, so that 
$2\beta > \frac{4}{3}$ in the required decay condition.
 
To the best knowledge of the authors, there are no publications
for the uniqueness $u=0$ by means of condition \eqref{1.7} requiring the 
exponential decay only in the limited coordinate components 
$(x_1, ..., x_m)$ and the polynomial decay 
in the rest components $(x_{m+1}, ..., x_n)$ 
with any choice of $m\in \{ 1, 2, ..., n\}$.
The key is a Carleman estimate (Lemma 2.1) in Section 2
where the weight function depends only on $x_1, ..., x_m$.

In Theorem \ref{T1}, the assumption means a fast exponential decay of
the solution $u(x)$ in $m$ components of the variable 
$x$ and a polynomial decay 
in the rest $(n-m)$ components of $x$.
Moreover, such critical exponential decay 
rate $2\beta$ is smaller if the exponent $\gamma$ of the behavior of $V(x)$ 
is larger, but $2\beta$ cannot be smaller or equal to $1$.
For larger $m$, the polynomial decay condition in $x_{m+1}, ..., x_n$ becomes 
more generous.
In particular, in the case of $m=1$, 
we can restrict the decay condition \eqref{1.7} to
one side $x_1 > 1$.  More precisely, we can prove
\\
\begin{thm}\label{T2}
{\it
Let $g$ satisfy \eqref{1.6} with $m=1$, and $\alpha$ be given by
(1.5).  We choose 
$r_1>0$ such that $\{ x\in \R^n;\, x_1 > r_1\} \subset \UUUUU$.
We assume that $u \in H^2_{loc}(\UUUUU)$ satisfies \eqref{1.4} 
and there exist constants 
$\beta > \alpha$ and $C>0$, $C_1>0$ such that  
\begin{equation}
\vert u(x)\vert \le Cg((x_2^2 + \cdots + x_n^2)^{\frac{1}{2}}
\exp( -C_1 \vert x_1\vert^{2\beta}), \quad \mbox{$(x_2,..., x_n) \in \R^{n-1}$
and $x_1 > r_1$ or $x_1 < -r_1$}.          \label{1.9}        
\end{equation}
Then $u =0$ in $\UUUUU$.
}
\end{thm}
Note that we assume \eqref{1.9} for $x_1 > r_1$ or $x_1 < -r_1$, but not for 
$\vert x_1\vert > r_1$.

Roughly speaking,
the existing results (e.g., \cite{KSW}, \cite{LMNN}) in the 
two-dimensional case $n=2$ are more refined than the case $n\ge 3$.
However, our main technique by a Carleman estimate can treat any dimensions
$n$, but, restricted to the case of $n=m=2$, our result is not the best 
possible.
\\

We emphasize that replacing an inequality \eqref{1.4} by an elliptic 
equation \eqref{1.1}, we can 
immediately obtain the same conclusions as Theorems \ref{T1} and \ref{T2}, 
Corollaries \ref{C1} and 1.3.
%\ref{C2}.  
More precisely, in \eqref{1.1} we consider that there exists a constant $C>0$
such that  
\begin{equation}
\vert A(x)\vert \le C\YYYY^{-\theta}, \quad 
\vert V(x)\vert \le C\YYYY^{-\gamma}, \quad x\in \UUUUU
                                                    \label{1.10}
\end{equation}
with constants $\gamma, \theta \in \R$.  
We recall that the constant $\alpha > 0$ satisfies 
\eqref{1.5}.  Then,
\\
\begin{thm}\label{T3}
{\it
If $u\in H^2_{loc}(\UUUUU)$ satisfies (1.1) with (1.10) and \eqref{1.7} with 
$\beta > \alpha$, then $u=0$ in $\UUUUU$.
}
\end{thm}

\begin{cor}\label{C2}

{\it
Let $m=n$ and let $u \in H^2_{loc}(\UUUUU)$ satisfy \eqref{1.1} with (1.10) and
\eqref{1.8} with $\beta > \alpha$.  Then $u=0$ in $\UUUUU$.
}
\end{cor}

\begin{cor}\label{C3}
{\it
Let $m \in \{ 1, ..., n\}$ be fixed and let $A=0$ in $\UUUUU$ and 
$V \in L^{\infty}(\UUUUU)$.
We assume that $u \in H^2_{loc}(\UUUUU)$ satisfies 
$$
\Delta u(x) - V(x)u(x) = 0 \quad \mbox{in $\UUUUU$}. 
$$
If there exist constants $C, C_1 > 0$ such 
that
$$
\vert u(x)\vert \le Cg(\vert z\vert)\exp(-C_1\vert y\vert^{\frac{4}{3}+\ep}),
\quad x \in \UUUUU \quad \mbox{with some $\ep > 0$},     
$$
then $u=0$ in $\UUUUU$.
}
\end{cor}

The case $m=n$ in Corollary 1.7 corresponds to the uniqueness result 
by Meshkov \cite{M} under (1.3).
 
\begin{thm}\label{T4}
{\it
We assume that $u \in H^2_{loc}(\UUUUU)$ satisfy \eqref{1.1} with (1.10). 
Let $m=1$ and let $\{ x\in \R^n;\, x_1 > r_1\} \subset \UUUUU$
with some $r_1 >0$.
Then \eqref{1.9} with $\beta > \alpha$, implies $u=0$ in $\UUUUU$.
}
\end{thm}

The article is composed of five sections.  Sections \ref{S2} and \ref{S3} are devoted to the 
proofs of Theorems \ref{T1} and \ref{T2} respectively.  In Section \ref{S4}, we provide that
special but not trivial choices $A$ and $V$ can imply $u=0$ in $\UUUUU$ if 
$u$ decays only with polynomial rate.  Section \ref{S5} provides the proof of
the key Carleman estimate: Lemma 2.1.

\section{Proof of Theorem \ref{T1}.}\label{S2}

We recall that $x=(x_1,..., x_n) = (y,z)$, where
$y\in \R^m$ and $z\in \R^{n-m}$.

Since $U \subset \R^n$ is a bounded domain, we can choose sufficiently large
$R_0>0$ such that $\ooo{U} \subset \{ x\in \R^n;\, \vert x\vert < R_0\}$, 
and so $\{ x\in \R^n;\, \vert y\vert > R_0\} \subset \R^n \setminus \ooo{U}$.
Hence, since $\{(y,z)\in \R^m \times \R^{n-m};\, \vert y\vert > R_0\}
\subset \{ x\in \R^n;\, \vert x\vert > R_0\}$, we see that 
$\{ y\in \R^m;\, \vert y\vert > R_0\} \times \R^{n-m} 
\subset \R^n \setminus \ooo{U}$. Moreover, without loss of generality, 
we can assume that $R_0=1$. 
If we finish the proof that $u\equiv 0$ in $\{ \vert y\vert > 1\} 
\times \R^{n-m}$, then the unique continuation (e.g., Choulli \cite{C})
yields that $u = 0$ in each component of $\UUUUU$ which is connected with 
the point at infinity in $\R^n$.  Since $U$ is simply connected, such a 
connected component is nothing but $\UUUUU$.
Therefore, $u \equiv 0$ in $\UUUUU$.  Thus it suffices to argue 
in the domain $\RRRMN$, where $y:= (x_1, ..., x_m)$, 
$z:= (x_{m+1}, ..., x_{n})$ and $x:= (y,z)$.
We recall that $\alpha$ is defined by \eqref{1.5} and $\alpha < \beta$
is satisfied.
\\
{\bf First Step.}
\\
For the constant $\alpha > 0$ defined by (1.5), we set 
\begin{equation}\label{function}
\va(x) = \vert y\vert^{2\alpha}, \quad x:=(y,z) 
\in \RRRMN.                    
\end{equation}
Then,
\begin{lem} \label{P1} (Carleman estimate).
{\it
%Let $u \in H^2_{loc}(\RRRMN)$ satisfy 
%\begin{equation*}
%\vert \Delta u(x) \vert \le C_0(\YYYY^{-\gamma} \vert u(x)\vert 
%+ \YYYY^{-\theta} \vert \nabla u(x) \vert + \vert F(x)\vert),
%\quad x\in \UUUUU,
%\end{equation*}
%where $F \in L^2_{loc}(\RRRMN)$.
Let $D \subset \RRRMN$ be an arbitrarily chosen bounded domain with 
piecewise smooth boundary $\ppp D$.  Then, there exist constants 
$s_0 > 0$ and $C>0$ such that 
\begin{equation*}
\int_D (s^3\vert y\vert^{6\alpha-4}\vert u\vert^2 
+ s\vert y\vert^{2\alpha-2} \vert \nabla u\vert^2) \weight dx
\le C\int_D \vert \Delta u\vert^2 \weight dx
\end{equation*}
for all $s\ge s_0$ and $u \in H^2_0(D)$.
}
\end{lem}

Here and henceforth we remark that $C>0$ and $s_0 > 0$ denote
generic constants which are 
independent of the choices of piecewise smooth bounded domains 
$D \subset \RRRMN$, and
are dependent on $C_0>0$ and $\theta, \gamma \in \R$.
\\
The proof of Lemma 2.1 is standard as elliptic Carleman estimate, but 
for demonstrating the uniformity of the estimate with respect to $D$, 
in Section 5, we provide a
direct proof of Lemma 2.1.  As for the recipe of the proof,
we can refer for example to Chapter 8 in Yamamoto \cite{Ya}, and see also 
Choulli \cite{C}.

Henceforth we set 
\begin{equation}
\sigma_1:= \max\{ 2\vert \theta\vert, \, \, 2\vert \gamma\vert,
\,\, 4\alpha\}.              \label{2.1}         
\end{equation}
Next, we will prove 
\begin{lem}\label{L1}
{\it
Let $D \subset \RRRMN$ be a bounded domain with piecewise smooth 
$\ppp D$, and let $q \in C^2_0(D)$ and $q \ge 0$ in $D$.
We assume that $u\in H^2(D)$ satisfies \eqref{1.4} in $\RRRMN$.  Then
\begin{equation*}
\int_D q(x) \vert \nabla u(x)\vert^2 dx 
\le C\int_D (\vert \Delta q(x)\vert 
+ \YYYY^{\sigma_1} q(x)) \vert u(x)\vert^2 dx.
\end{equation*}
Here $C>0$ is independent of choices of $q$ and piecewise smooth bounded 
domains $D \subset \RRRMN$.
}\end{lem}

Suitable choices of $q$ produce necessary inequalities for the proof of 
Theorem \ref{T1}.
\\
{\bf Proof of Lemma \ref{L1}.}
\\
By $q \in H^2_0(D)$, the integration by parts yields
\begin{equation*}
-\int_D uq\Delta u dx = \int_D q\nabla u\cdot \nabla u 
+ \int_D u\nabla q \cdot \nabla u dx,
\end{equation*}
and so
\begin{align*}
& \int_D q\vert \nabla u\vert^2 dx 
= - \frac{1}{2}\int_D \nabla q \cdot \nabla (u^2) dx 
- \int_D qu\Delta u dx\\
=& \frac{1}{2}\int_D (\Delta q) u^2 dx 
- \int_ D uq\Delta u dx.
\end{align*}
By \eqref{1.4}, we have 
\begin{align*}
& \int_D q\vert \nabla u\vert^2 dx 
\le \frac{1}{2}\int_D \vert \Delta q\vert u^2 dx
+ \int_D q\vert u\vert \vert \Delta u\vert dx\\
\le& \frac{1}{2}\int_D \vert \Delta q\vert u^2 dx
+ C_0\int_D q\vert u\vert (\YYYY^{-\gamma}\vert u\vert
+ \YYYY^{-\theta}\vert \nabla u\vert) dx.
\end{align*}
Moreover, for small $\ep > 0$, we can find a constant $C_{\ep} > 0$
such that 
\begin{align*}
& \int_D q\vert u\vert \YYYY^{-\theta}\vert \nabla u\vert dx
= \int_D \sqrt{q}\vert \nabla u\vert \sqrt{q} \YYYY^{-\theta}\vert u\vert dx
                                                       \\
\le& \ep\int_D q \vert\nabla u\vert^2 dx 
+ C_{\ep}\int_D q\YYYY^{-2\theta}\vert u\vert^2 dx.
\end{align*}
Hence,
\begin{align*}
& \int_D q\vert \nabla u\vert^2 dx 
\le \frac{1}{2}\int_D \vert \Delta q\vert u^2 dx
+ C_0\int_D q \YYYY^{-\gamma} \vert u\vert^2 dx\\
+ & C_0\ep\int_D q\vert \nabla u\vert^2 dx 
+ C_0C_{\ep}\int_D q\YYYY^{-2\theta}\vert u\vert^2 dx.
\end{align*}
Choosing $\ep > 0$ sufficiently small, we can absorb the third term 
on the right-hand side into the left-hand side, so that by $\YYYY \ge 1$
in $D$, we complete the proof of Lemma \ref{L1}.
$\blacksquare$
\\
{\bf Second Step.}
\\
For the proof, we need cut-off functions defined as follows.
We arbitrarily choose $\chi_0 \in C^2(\R)$ such that $0 \le \chi_0 \le 1$ and
$$
\chi_0(\xi) = 
\left\{ \begin{array}{rl}
&1, \quad \xi \ge 1, \\
&0, \quad \xi < 0.
\end{array}\right.
$$
We note that $\chi_0'(\xi) \ne 0$ only if $0\le \xi \le 1$.
Then, for constants $\rho_1, \rho_2$ satisfying 
$1 \le \rho_1 < \rho_1+1 < \rho_2$, we define
$$
\chi_{\rho_1+1,\rho_2}(\xi) : = \chi_0(\xi-\rho_1) - \chi_0(\xi - \rho_2).
$$
Then 
\begin{equation}
\chi_{\rho_1+1,\rho_2}(\xi) = 
\left\{ \begin{array}{rl}
&0, \quad \xi \le \rho_1, \\
&1, \quad \rho_1+1 \le \xi \le \rho_2, \\
&0, \quad \xi \ge \rho_2+1
\end{array}\right.
                                   \label{2.2} 
\end{equation}
and 
$\chi_{\rho_1+1,\rho_2} \in C^2(\R),$ $\mbox{supp}\,\chi_{\rho_1+1,\rho_2} \subset [\rho_1,\, \rho_2+1].$

Next, for a constant $\rho_3>0$, we set 
$\mu_{\rho_3}(\xi):= 1 - \chi_0(\xi-\rho_3)$.
Then, $\mu_{\rho_3} \in C^2(\R)$ and
\begin{equation}
\mu_{\rho_3}(\xi) = 
\left\{ \begin{array}{rl}
&1, \quad \xi \le \rho_3, \\
&0, \quad \xi \ge \rho_3+1. 
\end{array}\right.                     
                              \label{2.3}  
\end{equation}
For $n_1 \in \{1,2, ..., n\}$, we set $\xi:= (\xi_1, ..., \xi_{n_1})
\in \R^{n_1}$ and 
$\eta := (\eta_1, ..., \eta_{n-n_1}) \in \R^{n-n_1}$.
For $n= n_1$, we understand that $\xi = (\xi_1, ..., \xi_n) \in \R^n$ 
and the components $\eta$ do not appear. 
Then,
\begin{lem}\label{L2}
{\it
Let $1\le \rho_1 < \rho_1+1 < \rho_2$,
$1\le \rho_3 < \rho_3+1 < \rho_4$, and $r_1, r_2 > 0$.
\\
(i) We have
\begin{equation}
\mu_{r_1}(\vert \xi\vert)\mu_{r_2}(\vert \eta\vert) = 1
\quad \mbox{for $\vert \xi\vert \le r_1$ and $\vert \eta \vert \le r_2$}
                                       \label{2.4}   
\end{equation}
and
\begin{equation}
\sup_{r_1, r_2 > 0}
\Vert \ppp_i(\mu_{r_1}(\vert \xi\vert)\mu_{r_2}(\vert \eta\vert))
\Vert_{W^{1,\infty}(\R^n)} < \infty \quad \mbox{for $i \in \{1,.., n\}$}.
                                         \label{2.5}                  
\end{equation}
\\
(ii) We have
\begin{equation}
\chi_{\rho_1+1,\rho_2}(\vert \xi\vert)\mu_{r_1}(\vert \eta\vert) = 1
\quad \mbox{for $\rho_1+1 \le \vert \xi\vert \le \rho_2$ 
and $\vert \eta \vert \le r_1$}          \label{2.6}                    
\end{equation}
and
\begin{equation}
\sup_{1 \le \rho_1, \rho_1+1<\rho_2, r_1>0}
\Vert \ppp_i (\chi_{\rho_1+1,\rho_2}(\vert \xi\vert) \mu_{r_1}
(\vert \eta\vert))\Vert_{W^{1,\infty}(\R^n)} < \infty \quad 
\mbox{for $i \in \{1, ..., n\}$}.               
                                          \label{2.7}                   
\end{equation}
\\
(iii) We have
\begin{equation}
\chi_{\rho_1+1,\rho_2}(\vert \xi\vert)\chi_{\rho_3+1,\rho_4}(\vert \eta\vert) 
= 1 \quad \mbox{for $\rho_1+1 \le \vert \xi\vert \le \rho_2$ 
and $\rho_3+1 \le \vert\eta \vert \le \rho_4$}            \label{2.8}   
\end{equation}
and
\begin{equation}
\sup_{1 \le \rho_1, \rho_1+1<\rho_2, 1\le \rho_3, \rho_3+1<\rho_4}
\Vert \ppp_i (\chi_{\rho_1+1,\rho_2}(\vert \xi\vert) 
\chi_{\rho_3+1,\rho_4}(\vert \eta\vert))\Vert_{W^{1,\infty}(\R^n)} < \infty  
                                                         \label{2.9}     
\end{equation}
for $i=1,2, ..., n$.
}
\end{lem}
{\bf Proof of Lemma \ref{L2}.}
\\
Here and henceforth, we write $\chi'(\xi) = \frac{d\chi}{d\xi}(\xi)$ and
$\chi''(\xi) = \frac{d^2\chi}{d\xi^2}(\xi)$.
\\
By the definitions (2.3) and (2.4), 
we readily see \eqref{2.4}, \eqref{2.6} and \eqref{2.8}.
As for the uniform boundedness \eqref{2.5}, we can verify for 
example by 
\begin{align*}
& \frac{\ppp}{\eta_k}\mu_{r_1}(\vert \eta\vert) 
= \chi_0'(\vert \eta\vert - r_1)
\frac{\eta_k}{\vert \eta\vert}, \\
& \frac{\ppp^2}{\ppp \eta_j\ppp \eta_k}\mu_{r_1}(\vert \eta\vert) 
= \chi_0''(\vert \eta\vert-r_1)
\frac{\eta_j\eta_k}{\vert \eta\vert^2} 
+ \chi_0'(\vert \eta\vert - r_1)\frac{\vert \eta\vert^2\delta_{jk}
- \eta_j\eta_k}{\vert \eta\vert^3}
\end{align*}
for $1\le j, k \le n$.  Here $\delta_{jk} = 0$ if $j\ne k$ and 
$\delta_{kk} = 1$.
Therefore, since all the derivatives in the representations of 
the derivatives $\frac{\ppp}{\eta_k}\mu_{r_1}(\vert \eta\vert)$ and
$\frac{\ppp^2}{\ppp \eta_j\ppp \eta_k}\mu_{r_1}(\vert \eta\vert)$ 
vanish at $\vert \eta\vert =0$
by the definition of $\mu_{r_1}$.
The proofs of (2.8) and (2.10) are similar in terms of (2.3).
Thus we can complete the proof of Lemma \ref{L2}.
$\blacksquare$

We recall that $\va$ is defined by (2.1).  Moreover, for later uses, we prove
\begin{lem}\label{L3}
{\it
Let $q(x) = \chi(x)e^{2s\va}$, where 
\begin{equation*}
\chi(x):= \chi_{\rho_1,\rho_2}(\YYYY)\chi_{\rho_3,\rho_4}(\ZZZZ) 
\YYYY^{\sigma_1} \quad
\mbox{or} \quad  \chi_{\rho_1,\rho_2}(\YYYY)\mu_{r_1}(\ZZZZ)
\YYYY^{\sigma_1},
\end{equation*}
and $\sigma_1 > 0$ is defined by \eqref{2.1}, $1\le \rho_1 < \rho_1+1 
\le \rho_2$, $1\le \rho_3 < \rho_3+1 \le \rho_4$ and $r_1 >0$.
Then, 
\begin{equation*}
\vert \Delta q(x)\vert + 
\vert q(x)\vert \le Cs^2\YYYY^{4\alpha+\sigma_1}e^{2s\va(x)} \quad 
\mbox{for $x \in \RRRMN$ and $s\ge 1$.}
\end{equation*}
Here $C>0$ is independent of $s$, $x$ and $\sigma_1 > 0$ and
$\rho_1,\rho_2,\rho_3,\rho_4, r_1$
satisfying $1\le \rho_1 < \rho_1+1 \le \rho_2$,
$1\le \rho_3 < \rho_3+1 \le \rho_4$ and $r_1 >0$.}
\end{lem}

{\bf Proof of Lemma \ref{L3}.}
\\
Noting that $x\in \RRRMN$ implies $\YYYY > 1$, we directly see that 
\begin{equation*}
\vert \Delta \va(x)\vert+\vert \nabla \va(x)\vert +
\vert \va(x)\vert \le C\YYYY^{2\alpha} \quad 
\mbox{for $x \in \RRRMN$.}
\end{equation*}
Therefore, for $1\le i \le m$, we have 
\begin{equation*}
\vert \ppp_i(e^{2s\va(x)})\vert 
= \vert 2s(\ppp_i\va)e^{2s\va}\vert \le Cs\YYYY^{2\alpha}e^{2s\va}
\end{equation*}
and
\begin{equation}
 \vert \ppp_i^2(e^{2s\va(x)})\vert 
= \vert (4s^2(\ppp_i\va)^2 + 2s(\ppp_i^2\va)) e^{2s\va}\vert 
\le C(s^2\YYYY^{4\alpha} + s\YYYY^{2\alpha-2})e^{2s\va}
\le Cs^2\YYYY^{4\alpha}e^{2s\va}.
                                         \label{2.10}   
\end{equation}
Here for $m+1 \le i \le n$, we note that $\ppp_i(\weight) = 0$.
For short descriptions, we write
\begin{equation*}
\psi(\YYYY, \ZZZZ) := \chi_{\rho_1,\rho_2}(\YYYY)\chi_{\rho_3,\rho_4}
(\ZZZZ)
\end{equation*}
or
\begin{equation*}
\psi(\YYYY, \ZZZZ) := \chi_{\rho_1,\rho_2}(\YYYY)\mu_{r_1}(\ZZZZ).
\end{equation*}
Then, for $1 \le i \le m$, noting that $y=(x_1, ..., x_m)$,
we can directly calculate
\begin{equation*}
\ppp_i(\YYYY) = \frac{x_i}{\vert y\vert}, \quad 
\ppp_i^2(\YYYY) = \frac{\YYYY^2 - x_i^2}{\YYYY^3},
\end{equation*}
and so
\begin{align*}
& \ppp_i(\YYYY^{\sigma_1}) = \sigma_1\YYYY^{\sigma_1-2}x_i, \\
& \ppp_i^2(\YYYY^{\sigma_1}) = \sigma_1(\sigma_1-2)
\YYYY^{\sigma_1-4}x_i^2 + \sigma_1\YYYY^{\sigma_1-2}.
\end{align*}
Consequently, using also $\YYYY \ge 1$ and applying Lemma 2.3, we can 
obtain
\begin{equation}
\vert \ppp_i(\psi(\YYYY,\ZZZZ)\YYYY^{\sigma_1})\vert, \quad
\vert \ppp_i^2(\psi(\YYYY,\ZZZZ)\YYYY^{\sigma_1})\vert 
\le C\YYYY^{\sigma_1}                     \label{2.11}    
\end{equation}
for $1\le i \le m$.  
For $m+1 \le i \le n$, we can estimate similarly.
Since $q(x) = \psi(\YYYY,\ZZZZ)\YYYY^{\sigma_1} \weight$,
estimates \eqref{2.10} and \eqref{2.11} complete the proof of Lemma \ref{L3}.
$\blacksquare$
\\
{\bf Third Step.}
\\
For each $R > 6$, we consider 
\begin{equation*}
D(R):= \{ 4 \le \YYYY \le R\} \times \{ \ZZZZ \le R\}
= \{ x=(y,z) \in \R^m \times \R^{n-m};\, 
4 \le \YYYY \le R, \, \, \ZZZZ \le R\}.
\end{equation*}
We set 
\begin{align*}
&D_1(R):= \{ R-1 \le \YYYY \le R\} \times \{ \ZZZZ \le R\},\quad
D_2(R):= \{ 4 \le \YYYY \le R-1\} \times \{ R-1 \le \ZZZZ \le R\}, \\
& D_3(R):= \{ 4 \le \YYYY \le 5\} \times \{ \ZZZZ \le R-1\}, \quad
  D_0(R):= \{ 5 \le \YYYY \le R-1\} \times \{ \ZZZZ \le R-1\},
\end{align*}
$$
\vert \nabla( \chi_{5,R-1}(\YYYY)\mu_{R-1}(\ZZZZ)) \vert 
= \vert \Delta (\chi_{5,R-1}(\YYYY)\mu_{R-1}(\ZZZZ))\vert = 0
\quad \mbox{in $D_0(R)$}
$$
by (2.5), and 
\begin{equation*}
\www{u}(x):= \chi_{5,R-1}(\YYYY)\mu_{R-1}(\ZZZZ)u(x).
\end{equation*}
Then, 
$$
D(R) = D_1(R) \cup D_2(R) \cup D_3(R) \cup D_0(R)
$$
and
$\www{u} \in H^2_0(D(R))$ and 
\begin{equation*}
\www{u} = u \quad \mbox{in $\{ 5 \le \YYYY \le R-1\}
\times \{ \ZZZZ \le R-1\}$.}
\end{equation*}
Moreover \eqref{1.4} implies
\begin{align*}
& \Delta \www{u} 
= \chi_{5,R-1}(\YYYY)\mu_{R-1}(\ZZZZ)(\Delta u) \\
+ & 2\nabla(\chi_{5,R-1}(\YYYY)\mu_{R-1}(\ZZZZ))\cdot \nabla u
+   \Delta (\chi_{5,R-1}(\YYYY)\mu_{R-1}(\ZZZZ))u.
\end{align*}
Therefore, by (2.8) we obtain 
\begin{align*}
& \vert \Delta \www{u}(x)\vert \le C(\vert \Delta u(x)\vert 
+ \vert \nabla u(x)\vert + \vert u(x)\vert) \\
\le& C(1+\YYYY^{-\gamma})\vert u(x)\vert 
+ C(1+\YYYY^{-\theta})\vert \nabla u(x)\vert \\
\le& C\YYYY^{\frac{\sigma_1}{2}}\vert u(x)\vert 
+ C\YYYY^{\frac{\sigma_1}{2}}\vert \nabla u(x)\vert \quad
\mbox{in $D(R) \setminus D_0(R)$}.
\end{align*}
Here, we used that $\frac{\sigma_1}{2} 
\ge \max\{\vert \gamma\vert, \, \vert \theta\vert\}$ by \eqref{2.1}, and  
\begin{equation*}
(1+\YYYY^{-\gamma})^2 \le 2 + 2\YYYY^{-2\gamma}
\le C\YYYY^{\sigma_1},
\quad (1+\YYYY^{-\theta})^2 \le 2 + 2\YYYY^{-2\theta}
\le C\YYYY^{\sigma_1}   
\end{equation*}
for all $\YYYY \ge 1$.
Moreover, since (1.5) implies $6\alpha-4 \ge -2\gamma$ and 
$2\alpha -2 \ge -2\theta$, by $u = \www{u}$ in $D_0(R)$ and (1.4),
we have
\begin{align*}
& \vert \Delta \www{u}(x)\vert^2 
= \vert \Delta u(x)\vert^2 \\
\le& C_0(\YYYY^{-\gamma}\vert u(x)\vert 
+ \YYYY^{-\theta}\vert \nabla u(x)\vert)^2\\
\le& C(\YYYY^{6\alpha-4}\vert u(x)\vert^2 
+ \YYYY^{2\alpha-2}\vert \nabla u(x)\vert^2), \quad x\in D_0(R).
\end{align*}
Consequently, 
\begin{align*}
& \int_{D_0(R)} \vert \Delta \www{u}(x)\vert^2 e^{2s\va} dx \\
\le& C\int_{D_0(R)} (\YYYY^{6\alpha-4}\vert u(x)\vert^2 
+ \YYYY^{2\alpha-2}\vert \nabla u(x)\vert^2)e^{2s\va} dx\\
=&C\int_{D_0(R)} (\YYYY^{6\alpha-4}\vert \www{u}(x)\vert^2 
+ \YYYY^{2\alpha-2}\vert \nabla \www{u}(x)\vert^2)e^{2s\va} dx  \\
\le & C\int_{D(R)} (\YYYY^{6\alpha-4}\vert \www{u}(x)\vert^2 
+ \YYYY^{2\alpha-2}\vert \nabla \www{u}(x)\vert^2)e^{2s\va} dx.
\end{align*}
 
Hence, applying Lemma 2.1 to $\www{u}$ in 
$D(R)$, we obtain
\begin{align*}
&\int_{D(R)} (s^3\YYYY^{6\alpha-4}\vert \www{u}\vert^2
+ s\YYYY^{2\alpha-2}\vert \nabla \www{u}\vert^2) \weight dx
\le C\int_{D(R)} (\YYYY^{6\alpha-4}\vert \www{u}\vert^2 
+ \YYYY^{2\alpha-2}\vert \nabla \www{u}\vert^2)e^{2s\va} dx\\
+ & C\sum_{k=1}^3 \int_{D_k(R)} \YYYY^{\sigma_1} 
(\vert u\vert^2 + \vert \nabla u\vert^2) \weight dx.
\end{align*}
Choosing $s>0$ sufficiently large, we can absorb 
the first term on the right-hand side into the left-hand side, so that 
\begin{equation}   \label{2.12}           
 \int_{D(R)} (s^3\YYYY^{6\alpha-4}\vert \www{u}\vert^2
+ s\YYYY^{2\alpha-2}\vert \nabla \www{u}\vert^2) \weight dx
\le C\sum_{k=1}^3 \int_{D_k(R)} \YYYY^{\sigma_1} 
(\vert u\vert^2 + \vert \nabla u\vert^2) \weight dx
\end{equation}
$$
=: C(S_1(R) + S_2(R) + S_3(R)).
$$ 
\\
{\bf Estimation of $S_1$}
\\
For $D_1(R) = \{ R-1\le \YYYY \le R\} \times \{ \ZZZZ \le R\}$, we set
\begin{equation*}
\www{D}_1(R) := \{ R-2 \le \YYYY \le R+1\} \times \{ \ZZZZ \le R+1\}.
\end{equation*}
We set 
\begin{equation*}
q_1^R(x):= \chi_{R-1,R}(\YYYY)\mu_R(\ZZZZ)\YYYY^{\sigma_1}e^{2s\va} \ge 0.
\end{equation*}
Then, $q_1^R \in H^2_0(\www{D}_1(R))$.  Therefore, we can apply 
Lemmata 2.2 and 2.4 with $q_1^R$.  In terms of 
$q_1^R(x) = \YYYY^{\sigma_1}e^{2s\va(x)}$ for 
$x \in D_1(R)$, we have
\begin{align}
& \int_{D_1(R)} \YYYY^{\sigma_1}\vert \nabla u\vert^2 \weight dx
= \int_{D_1(R)} \vert \nabla u\vert^2 q_1^R dx 
\le \int_{\www{D}_1(R)} \vert \nabla u\vert^2 q_1^R dx \notag \\ 
\le & \int_{\www{D}_1(R)} ( \vert \Delta q_1^R\vert 
+ \YYYY^{\sigma_1}q_1^R) u^2 dx
\le C\int_{\www{D}_1(R)} s^2\YYYY^{\sigma_2} u^2 \weight dx.   \label{2.13} 
\end{align}
Here we set 
\begin{equation}
\sigma_2 := \sigma_1 + 4\alpha.                \label{2.XX}
\end{equation}
Here and henceforth $C>0$ denotes
generic constants which are independent of choices of $s \ge 1$ and $R > 6$.

Consequently, 
%by \eqref{1.8}, \eqref{1.9} 
by (1.6), (1.7) and $2\beta > 2\alpha$, for arbitrarily 
fixed $s>0$, applying \eqref{2.13}, we obtain
\begin{align*}
& \vert S_1(R)\vert 
= C\int_{D_1(R)} \YYYY^{\sigma_1}
(\vert u\vert^2 + \vert \nabla u\vert^2) \weight dx  \\
\le& C\int_{\www{D}_1(R)} s^2\YYYY^{\sigma_2} \vert u\vert^2 \weight dx\\
=& C\int_{\ZZZZ \le R+1} \int_{R-2 \le \YYYY\le R+1} s^2\YYYY^{\sigma_2}
g(\ZZZZ)^2 e^{-2C_1\YYYY^{2\beta}} e^{2s\YYYY^{2\alpha}} dydz \\
\le & C\left( \int_{\ZZZZ \le R+1} g(\ZZZZ)^2 dz\right) 
s^2(R+1)^{\sigma_2}e^{-2C_1(R-2)^{2\beta}}e^{2s(R+1)^{2\alpha}}
\int_{R-2 \le \YYYY \le R+1} dy \\
\le& Cs^2(R+1)^{\sigma_2}e^{-2C_1(R-2)^{2\beta}}e^{2s(R+1)^{2\alpha}}R^m
\, \longrightarrow \, 0 \quad \mbox{as $R \to \infty$}.
\end{align*}
For the first inequality, we used $\sigma_2 > \sigma_1$ and 
$\vert y\vert \ge 1$ for $(y,z) \in D_1(R)$, while for the 
last inequality we applied
$$
\int_{\vert z\vert\le R+1} g(\vert z\vert)^2 dx 
\le C\int^{R+1}_0 g(r)^2 r^{n-m-1} dr < \infty
$$
for all $R>0$ by means of (1.6).  
Hence, 
\begin{equation}
\lim_{R \to\infty} \vert S_1(R)\vert = 0 \quad \mbox{for arbitrarily fixed
$s>0$.}                  \label{2.14}        
\end{equation}
{\bf Estimation of $S_2$}
\\
For $D_2(R) = \{ 4 \le \YYYY \le R-1\} \times \{ R-1 \le \ZZZZ \le R\}$, we set
\begin{equation*}
\www{D}_2(R) := \{ 3 \le \YYYY \le R\} \times \{ R-2 \le \ZZZZ \le R+1\}
\end{equation*}
and
\begin{equation*}
q_2^R(x):= \chi_{4,R-1}(\YYYY)\chi_{R-1,R}(\ZZZZ)\YYYY^{\sigma_1}
e^{2s\va} \ge 0.
\end{equation*}
Then, $q_2^R \in H^2_0(\www{D}_2(R))$.

Similarly to \eqref{2.13}, Lemmata \ref{L1} and \ref{L3} yield
\begin{equation*}
\vert \Delta q_2^R(x)\vert +
\vert q_2^R(x)\vert \le Cs^2\YYYY^{\sigma_2}e^{2s\va(x)} \quad 
\mbox{for all $x\in \RRRMN$ and $s \ge 1$}.
\end{equation*}
Hence,
\begin{equation*}
 \int_{D_2(R)} \YYYY^{\sigma_1}\vert \nabla u\vert^2 \weight dx
\le C\int_{\www{D}_2(R)} s^2 \YYYY^{\sigma_2} u^2 \weight dx.
\end{equation*}
Therefore,
\begin{align*}
& \vert S_2(R)\vert 
= C\int_{D_2(R)} \YYYY^{\sigma_1}(\vert u\vert^2 + \vert \nabla u\vert^2) 
\weight dx 
\le C\int_{\www{D}_2(R)} s^2\YYYY^{\sigma_2} u^2\weight dx\\
\le& C\int_{R-2 \le \ZZZZ \le R+1} \int_{3\le \YYYY \le R} s^2\YYYY^{\sigma_2}
g(\ZZZZ)^2 e^{-2C_1\YYYY^{2\beta}} e^{2s\YYYY^{2\alpha}} dydz \\
\le & C\left( \int_{R-2 \le\ZZZZ \le R+1} g(\ZZZZ)^2 dz\right) 
\int_{\YYYY \ge 3} s^2\YYYY^{\sigma_2} e^{-2C_1\YYYY^{2\beta} 
+ 2s\YYYY^{2\alpha}}  dy\\
\le & C_2(s)\int^{R+1}_{R-2} g(r)^2 r^{n-m-1} dr.
\end{align*}
Here we used that 
\begin{equation*}
\int_{\YYYY\ge 1} s^2\YYYY^{\sigma_2} e^{-2C_1\YYYY^{2\beta} 
+ 2s\YYYY^{2\alpha}} dy =: C_2(s) < \infty
\end{equation*}
by $2\beta > 2\alpha$.
By \eqref{1.6}, we see 
$\lim_{R\to\infty} \int^{R+1}_{R-2} g(r)^2 r^{n-m-1} dr = 0$, so that 
\begin{equation}
\lim_{R \to\infty} \vert S_2(R)\vert = 0 \quad \mbox{for arbitrarily fixed
$s>0$.}               \label{2.15}           
\end{equation}
\\
{\bf Estimation of $S_3$}
\\
For $D_3(R) = \{ 4 \le \YYYY \le 5\} \times \{ \ZZZZ \le R-1\}$, we set
\begin{equation*}
\www{D}_3(R) := \{ 3 \le \YYYY \le 6\} \times \{ \ZZZZ \le R\}
\end{equation*}
and
\begin{equation*}
q_3^R(x):= \chi_{4,5}(\YYYY)\mu_{R-1}(\ZZZZ)\YYYY^{\sigma_1}e^{2s\va} \ge 0.
\end{equation*}
Then, using $\sigma_2 = \sigma_1 + 4\alpha$ by 
\eqref{2.XX} and $\vert y\vert \ge 1$ for $(y,z) \in \www{D}_3(R)$,
we obtain $q_3^R \in H^2_0(\www{D}_3(R))$ and 
$\vert \Delta q_3^R(x)\vert \le Cs^2\YYYY^{\sigma_2}\weight$.  
Similarly to \eqref{2.13},  
we apply Lemmata \ref{L1} and \ref{L3} with $q_3^R$ in $\www{D}_3(R)$:
\begin{align*}
& \int_{D_3(R)} \YYYY^{\sigma_1}\vert \nabla u\vert^2 \weight dx
= \int_{D_3(R)} \vert \nabla u\vert^2 q_3^R dx\\
\le& \int_{\www{D}_3(R)} \vert \nabla u\vert^2 q_3^R dx
\le C\int_{\www{D}_3(R)} s^2\YYYY^{\sigma_2} u^2\weight dx.
\end{align*}
Consequently, \eqref{1.6} and \eqref{1.7} yield 
\begin{align*}
& \vert S_3(R)\vert 
= C\int_{D_3(R)} \YYYY^{\sigma_1}(\vert u\vert^2 + \vert \nabla u\vert^2)
\weight dx
\le C\int_{\www{D}_3(R)} s^2\YYYY^{\sigma_2} u^2\weight dx\\
\le & \left( \int_{\ZZZZ \le R} g(\ZZZZ)^2 dz\right) 
\left( \int_{3\le \YYYY \le 6} s^2\YYYY^{\sigma_2} e^{-2C_1\YYYY^{2\beta} 
+ 2s\YYYY^{2\alpha}} dy \right)\\
\le& C\left(\int^R_0 g(r)^2 r^{n-m-1} dr \right)
 s^26^{\sigma_2} e^{-2C_13^{2\beta}}e^{2s 6^{2\alpha}} 
\left( \int_{3 \le \YYYY \le 6} dy \right)
\le Cs^2e^{2s 6^{2\alpha}},
\end{align*}
that is,
\begin{equation}
\vert S_3(R) \vert \le Cs^2e^{2s 6^{2\alpha}} \quad \mbox{for all large 
$s\ge 1$ and all $R>6$.}
                                        \label{2.16}    
\end{equation}
In view of \eqref{2.14} - \eqref{2.16}, estimate \eqref{2.12} yields
\begin{align*}
& \lim_{R\to\infty} \int_{5\le \YYYY \le R-1}
\int_{\ZZZZ \le R-1} \YYYY^{6\alpha-4}u^2 \weight dzdy
\le \lim_{R\to\infty} \int_{D(R-1)} \YYYY^{6\alpha-4}\vert \www{u}\vert^2 
\weight dx\\
\le& \lim_{R\to\infty} \int_{D(R)} \YYYY^{6\alpha-4}\vert \www{u}\vert^2 
\weight dx \le Cs^2e^{2s 6^{2\alpha}}
\end{align*}
for all large $s>0$.
Hence, 
\begin{equation*}
\int_{\YYYY \ge 5} \left(\int_{\R^{n-m}} \YYYY^{6\alpha-4}u^2 
e^{2s\YYYY^{2\alpha}} dz \right) dy \le Cs^2e^{2s 6^{2\alpha}}
\end{equation*}
for all large $s>0$.
Consequently, we obtain
\begin{equation*}
e^{2s 7^{2\alpha}} \int_{\YYYY \ge 7} \left( 
\int_{\R^{n-m}} \YYYY^{6\alpha-4}u^2 dz \right) dy  
\le \int_{\YYYY \ge 5} \left( \int_{\R^{n-m}} \YYYY^{6\alpha-4}u^2 
e^{2s\YYYY^{2\alpha}} dz\right) dy 
\le Cs^2e^{2s 6^{2\alpha}},
\end{equation*}
that is,
\begin{equation*}
\int_{\YYYY \ge 7} \int_{\R^{n-m}} \YYYY^{6\alpha-4}u^2 dzdy 
\le Cs^2e^{-2s (7^{2\alpha} - 6^{2\alpha})}
\end{equation*}
for all large $s>0$.  Letting $s \to \infty$, we see 
$u=0$ in $\{ \YYYY \ge 7\}\times \R^{n-m}$. 
The unique continuation yields $u=0$ in $\RRRMN$.
Thus the proof of Theorem \ref{T1} is complete.
$\blacksquare$

\section{Proof of Theorem \ref{T2}}\label{S3}

By $m=1$, we write $y:= x_1$ and $z:= (x_2, ..., x_n)$.
For each $R > 6$, we consider 
\begin{equation*}
D(R):= [4, R] \times \{ \ZZZZ \le R\}
= \{ x=(x_1,z) \in \R \times \R^{n-1};\, 
4 \le x_1 \le R, \, \, \ZZZZ \le R\}.
\end{equation*}
We set 
\begin{align*}
&D_1(R):= [R-1, R] \times \{ \ZZZZ \le R\},\quad
D_2(R):= [4,\, R-1] \times \{ R-1 \le \ZZZZ \le R\}, \\
& D_3(R):= [4,5] \times \{ \ZZZZ \le R-1\}.
\end{align*}
Then the same argument as in the proof of Theorem \ref{T1} can complete the
one of Theorem \ref{T2}.
$\blacksquare$
\\
{\bf Remark.}
\\
For Theorem \ref{T2}, the choice $m=1$ is essential.  For example, in the case of 
$m=2$, in order to establish a result similar to Theorem \ref{T2},
we have to consider a boundary neighborhood 
of $[4,R] \times [4,R] \times \{ \ZZZZ \le R\}$ and so estimate 
\begin{equation*}
\int_{\ZZZZ \le R-1}\left( \int^{R-1}_4 \left( \int_4^5 
(\vert u\vert^2 + \vert \nabla u\vert^2)\weight dx_1 \right) dx_2 
\right) dz.
\end{equation*}
However, we cannot control this integral by choosing large $R$, 
because of the integral over $[4,R-1]$ in $x_2$.
We note that such an integral does not appear in the case of $m=1$.

\section{Supplementary uniqueness results}\label{S4}
 
In this section, we consider elliptic equations.
Let $u \in H_{loc}^2(\R^n)$ satisfy 
\begin{equation}
\sumij \ppp_i(a_{ij}(x)\ppp_ju(x)) + A(x)\cdot \nabla u(x)
- V(x)u(x) = 0 \quad \mbox{in $\R^n$}.          \label{4.1}   
\end{equation}
Here we assume that 
\begin{equation}
\left\{ \begin{array}{rl}
& a_{ij} = a_{ji} \in C^1(\R^n)\quad \forall i,j\in\{1,\dots,n\},  \\
& \mbox{and there exists a constant $\kappa > 0$ such that}\\
& \sumij a_{ij}(x)\xi_i\xi_j \ge \kappa\sum_{i=1}^n \xi_i^2 \quad \mbox{for all 
$x\in \R^n$ and $\xi_1, ..., \xi_n \in \R$}.
\end{array}\right.
                               \label{4.2}          
\end{equation}
The main interest in the existing works
is the critical exponentially decay rate of $u$ 
which implies that $u=0$ in $\R^n$. 
As is referred in Section 1, this is related to the celebrated 
Landis' conjecture and there are many remarkable results.
Needless to say, the critical rate means that whenever {\bf each} $u$ decays
with such a rate, we can conclude that $u=0$ in $\R^n$. 
On the other hand, with special $A(x)$ and $V(x)$, we can prove
that much weaker rates imply that $u=0$ in $\R^n$ as follows.
\\
\begin{thm} \label{P2}
{\it 
Let $u\in H^2_{loc}(\R^n)$ satisfy \eqref{4.1} with (4.2).
Assume that $\ddd A \ge 0$ and $V \ge 0$ in $\R^n$ and $ A:= (a_1,..., a_n) 
\in C^1(\R^n),  V \in L^1_{loc}(\R^n)$.  
If there exists a sequence $R_\ell > 0$, $\ell\in \N$ such that   
$\lim_{\ell\to \infty} R_\ell = \infty$ and 
\begin{equation}
\lim_{\ell\to\infty} \frac{1}{R_\ell}
\int_{R_\ell \le \vert x\vert \le 2R_\ell} u^2 dx = 0,       \label{4.3}     
\end{equation}
then $u\equiv 0$ in $\R^n$.
}\end{thm}

\begin{cor}\label{C4}
{\it 
Assume that $\ddd A \ge 0$ and $V \ge 0$ in $\R^n$ and $ A:= (a_1,..., a_n)
\in C^1(\R^n),  V \in L^1_{loc}(\R^n)$. 
If a solution $u\in H^2_{loc}(\R^n)$ to \eqref{4.1} with (4.2) satisfies 
$$
\vert u(x)\vert \le C\vert x\vert^{-\frac{n}{2} + \ep} \quad
\mbox{for $\vert x\vert > 1$}
$$
with some $\ep < \frac{1}{2}$, then $u=0$ in $\R^n$.
}\end{cor}

For div $A \ge 0$ and $V\ge 0$ in $\R^n$, if we assume some decay 
conditions on both $u(x)$ and $\nabla u(x)$ as $\vert x\vert \to
\infty$, then multiplying (4.1) by $u$ and integrating by parts,
we can immediately prove $u=0$ in $\R^n$.  We emphasize that 
we never assume any decay conditions of $\vert \nabla u(x)\vert$. 

The proof is based on an energy estimate which is used for the 
Liouville problem for the Navier-Stokes equations
(e.g., the proof of Theorem X.9.5 on pp. 729-730
in Galdi \cite{G}).
\\
{\bf Proof of Theorem \ref{P2}.}
\\
{\bf First Step.}
\\
We choose $\psi_0 \in C^{\infty}(\R)$ such that $0\le \psi_0 \le 1$ and
\begin{equation*}
\psi_0(r) = 
\left\{ \begin{array}{rl}
& 1, \quad r \le 1, \\
& 0, \quad r\ge 2.
\end{array}\right.
\end{equation*}
For $R>0$, we set $\psi_R(r):= \psi_0\left( \frac{r}{R}\right)$.
Then, 
\begin{equation*}
\psi_R(r) = 
\left\{ \begin{array}{rl}
& 1, \quad r \le R, \\
& 0, \quad r\ge 2R.
\end{array}\right.
\end{equation*}
Moreover, we can prove: there exists a constant $C>0$, independent of $R>1$, 
such that
\begin{equation}
\left\{\begin{array}{rl}
& \ppp_i\ppp_j \psi_R(\vert x\vert) \le \frac{C}{R}\left(\frac{1}{R}
+ \frac{1}{\vert x\vert}\right), \quad x\in \R^n\setminus\{0\}, \\
& \Vert \nabla \psi_R(\vert x\vert)\Vert_{L^{\infty}(\R^n)}
\le \frac{C}{R} \quad \mbox{for $1\le i,j \le n$ and $R \ge 1$}.
\end{array}\right.
                                  \label{4.4}          
\end{equation}
\\
{\bf Verification of (4.4).}
\\
Direct calculations yield 
\begin{equation*}
\ppp_i\psi_R(\vert x\vert) = \frac{1}{R}\ppp_i(\vert x\vert)
\frac{d\psi_0}{dr}\left(\frac{\vert x\vert}{R}\right)
= \frac{1}{R} \frac{x_i}{\vert x\vert}\frac{d\psi_0}{dr}
\left(\frac{\vert x\vert}{R}\right)
\end{equation*}
and
\begin{equation*}
\ppp_i\ppp_j\psi_R(\vert x\vert) = \frac{1}{R^2} 
\frac{x_ix_j}{\vert x\vert^2}\frac{d^2\psi_0}{dr^2}
\left(\frac{\vert x\vert}{R}\right)
+ \frac{1}{R} \frac{\vert x\vert^2\delta_{ij} - x_ix_j}{\vert x\vert^3}
\frac{d\psi_0}{dr}\left(\frac{\vert x\vert}{R}\right)
\end{equation*}
for $1\le i,j \le n$, where we set $\delta_{ij} = 0$ if $i\ne j$ 
and $\delta_{ii} = 1$.
Therefore,
\begin{equation*}
\vert \ppp_i \psi_R(\vert x\vert) \vert \le \frac{C}{R}, \quad
\vert \ppp_i\ppp_j \psi_R(\vert x\vert)\vert \le \frac{C}{R}
\left( \frac{1}{R} + \frac{1}{\vert x\vert}\right)
\quad \mbox{for $x\in \R^n\setminus\{0\}$}.
\end{equation*}
The proof of \eqref{4.4} is complete.
$\blacksquare$
\\
{\bf Second Step.}
\\
We set $B_R:= \{ \vert x\vert < R\}$ and 
$\Sigma(R):= B_{2R} \setminus B_R$.
Multiply \eqref{4.1} in $B_{2R}$ with $u\psi_R$, and we have 
\begin{equation*}
-\int_{B_{2R}} \left(\sumij \ppp_i(a_{ij}\ppp_ju)\right) u\psi_R dx 
= \int_{B_{2R}} \sum_{k=1}^n a_k(\ppp_ku)u\psi_R dx 
 - \int_{B_{2R}} Vu^2\psi_R dx.
\end{equation*}
By $\psi_R = 0$ on $\ppp B_{2R}$, the integration by parts implies 
$$
 -\intBB \left( \sumij \ppp_i(a_{ij}\ppp_ju)\right) u\psi_R dx
= \intBB \sumij a_{ij}(\ppp_ju)(\ppp_iu)\psi_R dx 
+ \intBB \sumij a_{ij}(\ppp_ju)u\ppp_i\psi_R dx.
$$
%
%= & \frac{1}{2}\intBB \sumk a_k\ppp_k(u^2) \psi_R dx
%- \intBB Vu^2\psi_R dx.
%\end{align*}
Moreover, since the definition yields 
$\ppp_k\psi_R(\vert x\vert) = 0$ for $k=1,..., n$ and 
$\vert x\vert < R$ or $\vert x\vert > 2R$, $\psi_R(\vert x\vert) = 0$ for 
$\vert x\vert = 2R$, we see  
\begin{align*}
& \intBB \sumij a_{ij}(\ppp_ju)u\ppp_i\psi_R dx
= \frac{1}{2} \intBB \sumij a_{ij}\ppp_j(u^2)\ppp_i\psi_R dx\\
=& -\frac{1}{2}\intBB \sumij ((\ppp_ja_{ij})(\ppp_i\psi_R)
+ a_{ij}\ppp_i\ppp_j\psi_R)u^2 dx \\
=&  -\frac{1}{2} \int_{\Sigma(R)} \sumij ((\ppp_ja_{ij})(\ppp_i\psi_R)
+ a_{ij}\ppp_i\ppp_j\psi_R)u^2 dx,  
\end{align*}
so that (4.2) yields 
\begin{align*}
&  -\intBB \left( \sumij \ppp_i(a_{ij}\ppp_ju)\right) u\psi_R dx
\ge \intBB \kappa \vert \nabla u\vert^2 \psi_R dx \\
- & \frac{1}{2} \int_{\Sigma(R)} \sumij 
((\ppp_ja_{ij})(\ppp_i\psi_R) + a_{ij}(\ppp_i\ppp_j\psi_R)u^2 dx.
\end{align*}
Furthermore, again applying $\ppp_k\psi_R = 0$ on $\ppp B_{2R} \cup 
B_R$ and $\psi_R  = 0$ on $\ppp B_{2R}$ and 
$\ddd A \ge 0$ in $\R^n$, we obtain 
\begin{align*}
& \int_{B_{2R}} \sum_{k=1}^n a_k(\ppp_ku)u\psi_R dx 
= \frac{1}{2}\intBB \sumk a_k\ppp_k(u^2) \psi_R dx     \\
= & -\frac{1}{2}\intBB \sumk \ppp_k(a_k\psi_R)u^2 dx
+ \frac{1}{2} \int_{\ppp B_{2R}} \sumk a_k\psi_R u^2\nu_k dS \\
= & -\frac{1}{2}\intBB \left( \sumk \ppp_ka_k\right)
 \psi_Ru^2 dx
  -\frac{1}{2}\intBB \sumk a_k(\ppp_k\psi_R)u^2 dx
\le -\frac{1}{2}\int_{\Sigma(R)} \sumk a_k(\ppp_k\psi_R)u^2 dx
\end{align*}
Therefore, by $V \ge 0$ in $\R^n$, we reach 
\begin{align*}
& \kappa\intBB \vert \nabla u\vert^2 \psi_R dx 
\le \frac{1}{2}\int_{\Sigma(R)} u^2\left( \sumij (\ppp_ja_{ij})\ppp_i\psi_R
+ a_{ij}\ppp_i\ppp_j\psi_R)\right) dx\\
-& \frac{1}{2} \int_{\Sigma(R)} \sumk a_k(\ppp_k\psi_R)u^2 dx 
- \intBB Vu^2 \psi_R dx\\
\le& \frac{1}{2}\int_{\Sigma(R)} u^2\left( \sumij (\ppp_ja_{ij})\ppp_i\psi_R
+ a_{ij}\ppp_i\ppp_j\psi_R)\right) dx
- \frac{1}{2} \int_{\Sigma(R)} \sumk a_k(\ppp_k\psi_R)u^2 dx.
\end{align*}
Hence, since $\frac{1}{R} + \frac{1}{\vert x\vert} \le 2$ for 
$R > 1$ and $\vert x\vert \ge R$, estimate \eqref{4.4} implies  
\begin{align*}
& \kappa\intBB \vert \nabla u\vert^2 \psi_R dx 
\le \frac{1}{2}\int_{\Sigma(R)} u^2
\left\vert \left( \sumij (\ppp_ja_{ij})\ppp_i\psi_R
+ a_{ij}\ppp_i\ppp_j\psi_R)\right)\right\vert dx\\
+ & \frac{1}{2} \int_{\Sigma(R)} \left\vert \sumk a_k(\ppp_k\psi_R)u^2
\right\vert dx 
\le C\int_{\Sigma(R)} u^2\left( \frac{1}{R}
+ \frac{1}{R}\left( \frac{1}{R} + \frac{1}{\vert x\vert} \right)\right) dx
+ C\int_{\Sigma(R)} \frac{1}{R}u^2 dx \\
\le & \frac{C}{R}\int_{\Sigma(R)} u^2 dx.
\end{align*}
Hence, setting $R:= R_\ell$ for $\ell\in \N$ and using 
$\psi_{R_\ell} = 1$ in $B_{R_{\ell}}$, we have
\begin{equation*}
\int_{B_{R_\ell}} \vert \nabla u\vert^2 dx 
\le \frac{C}{R_\ell}\int_{\Sigma(R_\ell)}u^2 dx
= \frac{1}{R_\ell}o(R_\ell) = o(1)
\end{equation*}
as $\ell \to \infty$.  Consequently, $\nabla u = 0$ in $\R^n$, and so
$u$ is a constant function in $\R^n$. By \eqref{4.3}, we obtain
$u=0$ in $\R^n$.
$\blacksquare$
\\
{\bf Remark.}
In a similar estimation for the Navier-Stokes equations (e.g., \cite{G}), 
the extra condition
$\ddd A = 0$ is satisfied naturally 
by that the velocity field is a divergence free 
vector-valued function.
\\

In the case of $\Delta$, we can obtain sharper results than the general case
$\sumij \ppp_j(a_{ij}\ppp_i)$.
More precisely,  
\begin{thm}\label{P3}
{\it 
We assume that $u \in H_{loc}^2(\R^n)$ satisfies 
\begin{equation}
\Delta u - V(x)u(x) = 0 \quad \mbox{in $\R^n$},        \label{4.5}         
\end{equation}
and  
\begin{equation}
V \in L^1_{loc}(\R^n), \quad V\ge 0 \quad \mbox{in $\R^n$}.             
                              \label{4.6}    
\end{equation}
(i) If there exists a sequence $R_\ell$, $\ell\in \N$ such that 
$\lim_{\ell \to \infty}R_\ell = \infty$ and
\begin{equation*}
\lim_{\ell\to \infty} \frac{1}{R_\ell^2}
\int_{R_\ell\le \vert x\vert \le 2R_\ell} u^2 dx = 0,
\end{equation*}
then $u=0$ in $\R^n$.
\\
(ii) If there exists a sequence $R_\ell$, $\ell\in \N$ such that 
$\lim_{\ell \to \infty}R_\ell = \infty$ and
\begin{equation}
\lim_{\ell\to \infty} \int_{R_\ell-1\le \vert x\vert \le R_\ell+1} u^2 dx = 0,
                                             \label{4.7}
\end{equation}
then $u=0$ in $\R^n$.
}
\end{thm}

We note that the assumptions in (i) and (ii) are different
with the integral domains and the convergence orders of the integrals to 
$0$.
\begin{cor}\label{40}
{\it 
We assume that $u \in H^2_{loc}(\R^n)$ satisfies (4.5) with (4.6).  If 
\begin{equation*}
\vert u(x)\vert \le C\vert x\vert^{-\frac{n}{2}+\ep} \quad 
\mbox{for $\vert x\vert \ge 1$}
\end{equation*}
with some $\ep < 1$, then $u=0$ in $\R^n$.
}
\end{cor}

In Theorem \ref{P3} (ii), not assuming \eqref{4.7}, 
under the conditions $V \ge 0$ in $\R^n$ and 
$V \in L^1_{loc}(\R^n)$, we see that 
if $u \in H^2_{loc}(\R^n) \cap L^2(\R^n)$ satisfies (4.5), then 
$u\equiv 0$ in $\R^n$.  Indeed, (4.7) follows from $u \in L^2(\R^n)$.

We remark that the assumption in Corollary 4.4 does not necessarily imply 
$u\in L^2(\R^n)$, which means that $u\in L^2(\R^n)$ is not a necessary 
condition for 
$u\equiv 0$.  Moreover, for $0\le V(x)\le 1$, Corollary 4.4 improves
Theorem 1.1 in Das and Pinchover \cite{DP}.
\\
{\bf Proof of Theorem \ref{P3} (i)}.
\\
We recall that $B_R := \{ \vert x\vert < R\}$, and $\psi_R$ is defined in the 
proof of Theorem 4.1, and satisfies \eqref{4.4}.
We multiply \eqref{4.5} with $u\psi_R$, we have 
\begin{equation*}
-\int_{B_{2R}} (\Delta u)u\psi_R dx = -\int_{B_{2R}} Vu^2\psi_R dx.
\end{equation*}
Using $\psi_R(\vert x\vert) = \vert \nabla\psi_R(\vert x\vert)\vert 
= 0$ for $\vert x\vert = 2R$, by the integration by parts we obtain
\begin{align*}
& -\int_{B_{2R}} (\Delta u)u\psi_R dx 
= \int_{B_{2R}} \vert \nabla u\vert^2 \psi_R dx
 + \int_{B_{2R}} u\nabla u\cdot \nabla \psi_R dx \\
=& \int_{B_{2R}} \vert \nabla u\vert^2 \psi_R dx
 - \frac{1}{2}\int_{B_{2R}} u^2 \Delta \psi_R dx,
\end{align*}
which implies
\begin{align*}
& \int_{B_{2R}} \vert \nabla u\vert^2 \psi_R dx
= \frac{1}{2} \int_{B_{2R}} u^2 \Delta \psi_R dx
- \int_{B_{2R}} Vu^2\psi_R dx  \\
\le& \frac{1}{2} \int_{R\le \vert x\vert \le 2R} u^2 \Delta \psi_R dx
\end{align*}
by means of $\Delta \psi_R(x) = 0$ for $\vert x\vert \le R$ or
$\vert x\vert \ge 2R$, and $V \ge 0$ and $\psi_R \ge 0$ in 
$\R^n$.  Since $\psi_R(x) = 1$ in $B_R \subset B_{2R}$, in terms of 
\eqref{4.4}, we obtain
\begin{equation*}
\int_{B_R} \vert \nabla u\vert^2 dx \le \frac{C}{R^2}
\int_{R\le \vert x\vert \le 2R} u^2 dx.
\end{equation*}
Choosing $R:= R_\ell$ and letting $\ell\to \infty$, by the assumption, we have
$\int_{\R^n} \vert \nabla u\vert^2 dx = 0$, that is,
$u=c_0$ is a constant function.  Again the assumption implies $c_0=0$.
\\
{\bf Proof of Theorem \ref{P3} (ii).}
\\
Setting 
\begin{equation*}
\www{\psi}(r)=
\left\{ \begin{array}{rl}
&1, \quad r\le 1, \\
&0, \quad r\ge 3,
\end{array}\right.
\end{equation*}
we define a different cut-off function:
$\www{\psi}_R(r):= \www{\psi}(r-R+2)$.  Then,
\begin{equation*}
\www{\psi}_R(\vert x\vert) =
\left\{ \begin{array}{rl}
&1, \quad \vert x\vert < R-1, \\
&0, \quad \vert x\vert > R+1,
\end{array}\right.
\end{equation*}
and $\www{\psi}_R(\vert x\vert) = \vert \nabla\www{\psi}_R(\vert x\vert)\vert
= 0$ if $\vert x\vert = R+1$, and  
such that 
\begin{equation*}
\sup_{x\in \R^n, R>1} \vert \Delta (\www{\psi}_R(\vert x\vert))\vert <
\infty.
\end{equation*}
Indeed, 
\begin{equation*}
\ppp_i^2 \www{\psi}_R(x) =
\left\{ \begin{array}{rl}
&\www{\psi}''(\vert x\vert - R+2)\frac{x_i^2}{\vert x\vert^2}
+ \www{\psi}'(\vert x \vert - R+2)\frac{\vert x\vert^2 - x_i^2}
{\vert x\vert^3}, \\
&\qquad \quad \qquad \qquad \qquad \qquad \qquad \qquad 
 R-1 \le \vert x \vert \le R+1, \\
&0, \quad \vert x\vert \le R-1,\, \vert x\vert \ge R+1.
\end{array}\right.
\end{equation*}
Multiplying \eqref{4.5} with $u\www{\psi}_R$ and integrating by parts over 
$B_{R+1}$, similarly to the proof of Theorem \ref{P3} (i), 
in view of $V\ge 0$ in $\R^n$ and the conditions of $\www{\psi}_R$ on 
$\ppp B_{R+1}$ we obtain
\begin{equation*}
\int_{B_{R+1}} \www{\psi}_R \vert \nabla u\vert^2 dx 
\le \frac{1}{2} \int_{R-1\le \vert x\vert \le R+1} u^2\Delta \www{\psi}_R dx.
\end{equation*}
Since $\www{\psi}_R = 1$ for $\vert x\vert \le R-1$, we see
\begin{equation*}
\int_{B_{R-1}} \vert \nabla u\vert^2 dx 
\le C\int_{R-1\le \vert x\vert \le R+1} u^2 dx.
\end{equation*}
Setting $R=R_\ell$ and letting $\ell \to \infty$, we can reach the 
conclusion of (ii). 
$\blacksquare$

We can rewrite Theorem \ref{P3} (ii) as
\begin{thm}\label{P4}
{\it
Let $u \in H^2_{loc}(\R^n)$ satisfy \eqref{4.5} with \eqref{4.6} and 
$u \not\equiv 0$ in $\R^n$.  Then,
\begin{equation}
\inf_{R>1} \int_{R-1\le \vert x\vert \le R+1} u^2 dx > 0.
                                                         \label{4.8}   
\end{equation}
}
\end{thm}
{\bf Proof of Theorem \ref{P4}.}
\\
If the conclusion (4.8) does not hold, then we have two cases:
\\
{\bf (I):}
\begin{equation*}
\int_{R_0-1\le \vert x\vert \le R_0+1} u^2 dx = 0 \quad 
\mbox{with some $R_0 > 0$}.
\end{equation*}
{\bf (II):} We can find a sequence $R_\ell$, $\ell\in \N$ such that 
$\lim_{\ell\to\infty} R_\ell = \infty$ and 
\begin{equation*}
\lim_{\ell\to \infty} \int_{R_\ell-1\le \vert x\vert \le R_\ell+1} u^2 dx
= 0.
\end{equation*} 
In the former case, we obtain $u=0$ for $R_0-1\le \vert x\vert \le R_0+1$ and 
the unique continuation yields $u=0$ in $\R^n$. This is impossible.
In the latter case, Theorem 4.3 (ii) yields $u=0$ in $\R^n$, which is a 
contradiction.  Thus the proof is complete.
$\blacksquare$

In the same case, Theorem 1.2 in \cite{KSW} reads that there exists a 
constant $C>0$ such that 
\begin{equation*}
\inf_{\vert x_0\vert = R}\sup_{\vert x-x_0\vert<1} \vert u(x)\vert
\ge e^{-CR\log R} \quad \mbox{for all $R > 1$,}      
\end{equation*}
which is a different type of conclusion from \eqref{4.8}.

\section{Proof of Lemma \ref{P1}}\label{S5}

We recall that $\va(x) = \vert y\vert^{2\alpha}$ for 
$x := (y,z) \in \RRRMN$.
\\
{\bf First Step.}
We prove
\begin{lem}\label{L40}
{\it
We assume $\va \in C^2(\ooo{D})$ and $D \subset \RRRMN$ is 
a piecewise smooth bounded domain.  Then, there exist constants
$C>0$ and $s_0>1$ such that 
\begin{align*}
& \int_D \{ 4s^3\sumij (\ppp_i\ppp_j\va)(\ppp_i\va)(\ppp_j\va)
- s(\Delta^2\va)\} \vert u\vert^2 \weight dx
+ 4s\int_D \sumij (\ppp_i\ppp_j\va)\ppp_i(ue^{s\va})
\ppp_j(ue^{s\va}) dx\\
\le &C\int_D \vert \Delta u\vert^2 \weight dx 
\end{align*}
for all $s\ge s_0$ and $u\in H^2_0(D)$.
Here the constants $C$ and $s_0$ are independent of choices of piecewise 
smooth bounded domain $D \subset \RRRMN$. 
}\end{lem}
{\bf Proof of Lemma \ref{L40}.}
\\
For treating the weighted $L^2$-norms, we introduce
\begin{equation*}
w(x) = e^{s\va(x)}u(x), \quad
Pw(x) := e^{s\va}(- \Delta)(e^{-s\va}w).
\end{equation*}
Henceforth we write
\begin{equation*}
(u\cdot v):= \int_D u(x)v(x) dx, \quad 
\Vert v\Vert:= (v\cdot v)^{\frac{1}{2}} \quad \mbox{for $u,v \in L^2(D)$}.
\end{equation*}

Our goal is to obtain a lower estimate of $\Vert Pw\Vert^2$.
Direct calculations yield
\begin{equation*}
Pw = - \Delta w 
+ 2s\nabla\va \cdot \nabla w + (- s^2\vert \nabla \va\vert^2 + s\Delta \va)w.
                      \label{5.1}    
\end{equation*}
A traditional way for obtaining a lower estimate for 
$\Vert Pw\Vert^2$, is the decomposition of the operator
$P$ into the symmetric part $P_+$ and the antisymmetric part
$P_-$, that is, $Pw = P_+w + P_-w$.
We consider the formal adjoint operator $P^*$ to $P$:
\begin{equation*}
(Pu \cdot v) = (u \cdot P^*v), \quad 
u,v \in C^{\infty}_0(D).
\end{equation*}
Then, we can calculate as 
\begin{equation*}
P^*w = - \Delta w - 2s\nabla\va\cdot\nabla w
- (s\Delta \va + s^2\vert \nabla\va\vert^2)w.
\end{equation*}
We define the symmetric part $P_+$ and the antisymmetric part $P_-$
of $P$ by 
\begin{equation*}
P_+ = \frac{1}{2}(P+P^*), \quad
P_- = \frac{1}{2}(P-P^*).
\end{equation*}
Then, 
\begin{equation}
P_+w = -\Delta w - s^2\vert\nabla\va\vert^2 w, \quad
P_-w = 2s\nabla\va\cdot\nabla w + s(\Delta \va)w.   \label{5.1}
\end{equation}
Since $Pw = P_+w + P_-w$, we have 
\begin{equation} \label{5.2}    \begin{split}
&\int_D \vert \Delta u\vert^2 \weight dx = \Vert P_+w + P_-w\Vert^2
= \Vert P_+w\Vert^2 + \Vert P_-w\Vert^2\\
+ & 2(P_+w \cdot \, P_-w) \ge 2(P_+w \cdot\, P_-w) + \Vert P_+w\Vert^2.
  \end{split}                                           
\end{equation}
Therefore, we have
\begin{align*}
& \int_D \vert \Delta u\vert^2 e^{2s\va} dx 
 \ge 2(P_+w \cdot\, P_-w)
= 2((-\Delta w - s^2\vert\nabla\va\vert^2 w)\,\, \cdot\,\,
(2s(\nabla\va \cdot \nabla w) + s(\Delta \va)w))\\
=&  2(-\Delta w \cdot\, 2s(\nabla\va\cdot\, \nabla w)) 
+ 2(-\Delta w\cdot\, s(\Delta \va)w)
\end{align*}
\begin{equation}
- 2(s^2\vert\nabla\va\vert^2 w\cdot\, 2s(\nabla\va\cdot\, \nabla w))
- 2(s^2\vert\nabla\va\vert^2w\cdot\, s(\Delta \va)w)
=: \sum_{k=1}^4 J_k.     \label{5.3x}
\end{equation}
Henceforth $C>0$ denotes generic constants which are independent of $s$
and the choice of bounded domains $D \subset \RRRMN$, but depend on $\va$.
We note that $C$ may change line by line.

Applying $2(\ppp_kw)(\ppp_k\ppp_jw) = \ppp_j(\vert \ppp_kw\vert^2)$ 
and integration by parts, we have
\begin{align*}
& J_1 = -2\int_D (\Delta w) (2s\nabla\va \cdot \nabla w) dx
= - 4s\sum_{j,k=1}^n \int_D (\ppp_k^2w) (\ppp_jw)\ppp_j\va dx\\
= & 4s\sum_{j,k=1}^n \int_D (\ppp_kw) (\ppp_k\ppp_jw)(\ppp_j\va) dx 
+ 4s\int_D (\ppp_kw) (\ppp_jw)(\ppp_k\ppp_j\va) dx\\
=& -2s\int_D (\Delta \va)\vert \nabla w\vert^2dx
+ 4s\sum_{j,k=1}^n\int_D (\ppp_jw)(\ppp_kw)(\ppp_j\ppp_k\va) dx.
\end{align*}
Moreover, 
\begin{equation*}
 2s\int_D (\nabla (\Delta \va)\, \cdot w\nabla w) dx 
= s\int_D (\nabla (\Delta \va) \, \cdot\nabla (\vert w\vert^2)) dx
= -s\int_D (\Delta^2 \va)\vert w\vert^2 dx
\end{equation*}
Hence, the Green formula yields 
\begin{align*}
& J_2= -2\int_D (\Delta w) s(\Delta \va)w dx 
= 2s\int_D (\nabla w \cdot \nabla((\Delta\va)w)) dx  \\
=& 2s\int_D (\Delta\va)\vert \nabla w\vert^2 dx
+ 2s\int_D (\nabla(\Delta\va)\, \cdot w\nabla w) dx  \\
=& 2s\int_D (\Delta\va)\vert \nabla w\vert^2 dx
- s\int_D (\Delta^2\va)\vert w\vert^2 dx.
\end{align*}

Next, since
\begin{equation*}
\nabla (\vert \nabla \va\vert^2) \cdot \nabla \va 
= \sum_{i=1}^n \ppp_i\left(\sum_{j=1}^n \vert \ppp_j\va\vert^2\right)
\ppp_i\va
= \sum_{i,j=1}^n 2(\ppp_i\ppp_j\va)(\ppp_i\va)\ppp_j\va,
\end{equation*}
we obtain
\begin{align*} 
& J_3 = -2\int_D s^2\vert\nabla \va\vert^2 w \cdot 
2s(\nabla\va \cdot \nabla w) dx
= -4s\sum^n_{i=1} \int_D (s^2\vert\nabla \va\vert^2 w)
(\ppp_i\va)(\ppp_iw) dx\\
=& -2s\sum_{i=1}^n \int_D s^2\vert\nabla \va\vert^2 
(\ppp_i\va)\ppp_i(w^2) dx
= 2s \sum_{i=1}^n \int_D \ppp_i(s^2\vert\nabla \va\vert^2\ppp_i\va)w^2 dx\\
= & 2s\int_D \{ (\nabla(s^2\vert\nabla \va\vert^2) \cdot \, \nabla\va) 
+ s^2\vert\nabla \va\vert^2 \Delta\va\}w^2 dx\\
=&  4s^3\int_{D} \sumij (\ppp_i\ppp_j\va)(\ppp_i\va)(\ppp_j\va)
w^2 dx + 2s^3\int_D \vert \nabla\va\vert^2 (\Delta \va) w^2 dx. 
\end{align*}
Finally 
\begin{equation*}
 J_4 = -2s\int_D s^2\vert\nabla \va\vert^2w
(\Delta \va)w dx
= -2s^3 \int_D \vert \nabla \va\vert^2 (\Delta\va) w^2 dx. 
\end{equation*}
Therefore, 
\begin{align*}
& \sum_{k=1}^4 J_{k}
= 4s\int_D \sum_{i,j=1}^n (\ppp_i\ppp_j\va)(\ppp_iw)(\ppp_jw) dx
- s\int_D (\Delta^2\va) \vert w\vert^2 dx \\
+ &4s^3\int_D \sumij (\ppp_i\ppp_j\va)(\ppp_i\va)(\ppp_j\va)\vert w\vert^2 dx.
\end{align*}
Thus the proof of Lemma \ref{L40} is complete.
$\blacksquare$
\\
{\bf Second Step.}
\\
Let the function $\varphi$ be defined by (\ref{function}) and the constant 
$\alpha>\frac 12$ is defined by \eqref{1.5}.  
Then,
\begin{equation}
\left\{ \begin{array}{rl}
& \ppp_i\va = 2\alpha \YYYY^{2\alpha-2}y_i, \\
& \ppp_i\ppp_j\va = \left\{ \begin{array}{rl}
      &4\alpha(\alpha-1)\YYYY^{2\alpha-4}y_iy_j, \quad i\ne j,\\
      &4\alpha(\alpha-1)\YYYY^{2\alpha-4}y_i^2
        + 2\alpha\YYYY^{2\alpha-2}, \quad i=j,
     \end{array}\right.
\end{array}\right.
                                                       \label{5.3}        
\end{equation}
for $1\le i,j \le m$.
Therefore,
\begin{equation}\begin{split}\label{5.4}
& \sum_{i,j=1}^n (\ppp_i\ppp_j\va)(\ppp_i\va)(\ppp_j\va)   
 = \sum_{i,j=1}^m (\ppp_i\ppp_j\va)(\ppp_i\va)(\ppp_j\va)          
= \sum_{i,j=1}^m 4\alpha(\alpha-1)\YYYY^{2\alpha-4}y_iy_j 
4\alpha^2\YYYY^{4\alpha-4}y_iy_j  \\
+& 2\alpha \YYYY^{2\alpha-2}4\alpha^2\YYYY^{4\alpha-4}
\sum_{i=1}^m y_i^2
= 8\alpha^3(2\alpha-1)\YYYY^{6\alpha-4}. 
\end{split}\end{equation}
Moreover,
\begin{equation}
\vert \Delta^2\va(x)\vert \le C\YYYY^{2\alpha-4}, \quad 
x\in \RRRMN.                                       \label{5.5}        
\end{equation}
By (1.5), we have $\alpha > \frac{1}{2}$.

For all $s>0$, by (5.5) and (5.6), we see
$$
 4s^3 \sum_{i,j=1}^m (\ppp_i\ppp_j\va)(\ppp_i\va)(\ppp_j\va) 
- s\Delta^2\va \ge
$$
\begin{equation}\label{5.7X}
 32\alpha^3(2\alpha-1)s^3\YYYY^{6\alpha-4} - Cs\YYYY^{2\alpha-4}
\quad \mbox{in $\RRRMN$}.                   
\end{equation}
Henceforth, for $a:= (a_1,..., a_m)$, $b = (b_1, ..., b_m) \in \R^m$, we
write \\
$(a\cdot b) := \sum_{i=1}^m a_ib_i$ and 
$\vert a\vert = (a\cdot a)^{\frac{1}{2}}$.

Using $\alpha > \frac{1}{2}$ and the Cauchy-Schwartz inequality
$(a\cdot y)^2 \le \vert a\vert^2 \vert y\vert^2$, we have
\begin{align*}
& \sum_{i,j=1}^m (\ppp_i\ppp_j\va)a_ia_j
= 2\alpha \YYYY^{2\alpha-2}\vert a\vert^2
+ 4\alpha(\alpha-1)\YYYY^{2\alpha-4}\sum_{i,j=1}^m y_iy_ja_ia_j\\
=& 2\alpha\YYYY^{2\alpha-2}\vert a\vert^2
+ (2\alpha(2\alpha-1) - 2\alpha)\YYYY^{2\alpha-4}(a\cdot y)^2\ge \\
& 2\alpha\YYYY^{2\alpha-2}\vert a\vert^2
   - 2\alpha\YYYY^{2\alpha-4}(a\cdot y)^2
= 2\alpha\YYYY^{2\alpha-4}(\YYYY^2 \vert a\vert^2 - (a\cdot y)^2).
\end{align*}
Consequently, we have
\begin{equation}
4s\int_D \sum_{i,j=1}^m (\ppp_i\ppp_j\va)(\ppp_i(ue^{s\va}))
(\ppp_j(ue^{s\va})) dx 
\ge 0.                        \label{5.6}       
\end{equation}
Therefore, we choose $s_0>0$ sufficiently large such that
$$
32\alpha^3(2\alpha-1)s^3 \vert y\vert^{6\alpha-4}
- Cs\vert y\vert^{2\alpha-4}
\ge \alpha^3(2\alpha-1)s^3\vert y\vert^{6\alpha-4},\quad
\vert y\vert \ge 1 
$$
for $s \ge s_0$.  In terms of (5.7) and (5.8), Lemma 5.1 yields that 
there exist constants $s_0 \ge 1$ and $C>0$ such that 
\begin{equation}
\int_D  \alpha^3(2\alpha-1)s^3\YYYY^{6\alpha-4}\weight dx
\le C\int_D \vert \Delta u\vert^2 \weight dx   \label{5.7}    
\end{equation}
for all $s \ge s_0$ and all $u \in H^2_0(D)$.
\\
{\bf Third Step.}
\\
In addition to (5.9), we have to estimate also 
$\vert \nabla u\vert^2\weight$.
In view of (5.9), \eqref{5.2} and \eqref{5.3x}, we have 
\begin{equation}
\Vert P_+w\Vert^2
+ s^3\int_D \YYYY^{6\alpha-4}u^2 \weight dx
\le C\int_D \vert \Delta u\vert^2 \weight dx.
                                        \label{5.8}           
\end{equation}
Using (5.1) and (5.4), we establish 
\begin{align*}
& s\int_D (P_+w)\YYYY^{2\alpha-2}w dx
= s\int_D (-\Delta w)w \YYYY^{2\alpha-2} dx 
- s^3\int_D \vert \nabla \va\vert^2 \YYYY^{2\alpha-2} w^2 dx\\
=& s\int_D \vert \nabla w\vert^2 \YYYY^{2\alpha-2} dx 
+ s\int_D (w\nabla w\cdot y) 2(\alpha-1)\YYYY^{2\alpha-4} dx
- s^3\int_D 4\alpha^2 \YYYY^{6\alpha-4} w^2 dx.
\end{align*}
Moreover, using $\ppp_k(\YYYY^{2\alpha-4}) = (\alpha-2)\YYYY^{2\alpha-6}2y_k$
for $1\le k \le m$, we obtain
\begin{align*}
& s\int_D (w\nabla w\cdot y) 2(\alpha-1)\YYYY^{2\alpha-4} dx  \\
=& s(\alpha-1)\int_D \sum_{k=1}^m \ppp_k(w^2)y_k \YYYY^{2\alpha-4} dy\\
=& -s(\alpha-1)\int_D w^2 \sum_{k=1}^m \ppp_k(y_k \YYYY^{2\alpha-4}) dy\\
=& -s(\alpha-1)\int_D w^2(m\YYYY^{2\alpha-4} + 2(\alpha-2)\YYYY^{2\alpha-4})dy
\\
=& -s(\alpha-1)(m+2(\alpha-2))\int_D w^2\YYYY^{2\alpha-4} dy.
\end{align*}
Therefore, 
\begin{equation}\begin{split}\nonumber
& s\int_D \vert \nabla w\vert^2 \YYYY^{2\alpha-2} dx 
= s\int_D (P_+w)w\YYYY^{2\alpha-2} dx \\
+ &s(\alpha-1)(m+2(\alpha-2))\int_D w^2\YYYY^{2\alpha-4} dx
+s^3\int_D 4\alpha^2\YYYY^{6\alpha-4} w^2 dx.
\end{split}
\end{equation}
Applying \eqref{5.8} and 
\begin{equation*}
\left\vert s\int_D (P_+w)w \YYYY^{2\alpha-2} dx \right\vert
\le \frac{1}{2} \int_D \vert P_+w\vert^2 dx 
+ \frac{1}{2}s^2\int_D \YYYY^{4\alpha-4} w^2 dx,
\end{equation*}
we obtain
\begin{align}
& s\int_D \vert \nabla w\vert^2 \YYYY^{2\alpha-2} dx      
\le \frac{1}{2}\int_D \vert \Delta u\vert^2 \weight dx
   + \frac{1}{2} \int_D s^2\vert y\vert^{4\alpha-4}\vert w\vert^2 dx
                             \notag \\
+ & Cs\int_D  w^2\YYYY^{2\alpha-4} dx 
+ 4s^3\int_D 4\alpha^2 \YYYY^{6\alpha-4} w^2 dx.  \label{5.9} 
\end{align}

Finally, applying (5.9) and $(\nabla u)e^{s\va}
= -s(\nabla \va)w + \nabla w$ and 
$\nabla \va = 2\alpha\vert y\vert^{2\alpha-2}y$, 
we estimate the last term in (\ref{5.9}):

\begin{equation*}
\int_D (s^3\YYYY^{6\alpha-4}u^2 + s\YYYY^{2\alpha-2}
\vert \nabla u\vert^2) \weight dx
\le C\int_D \vert \Delta u\vert^2 \weight dx
\end{equation*}
for all $s \ge s_0$ and $u \in H^2_0(D)$.
Thus the proof of Lemma \ref{P1} is complete.
$\blacksquare$
\\

{\bf Acknowledgments.}
Masahiro Yamamoto was supported by 
Grant-in-Aid for Scientific Research (A) 20H00117 
and Grant-in-Aid for Challenging Research (Pioneering) 21K18142 of 
Japan Society for the Promotion of Science.

\end{document}